\documentclass[12pt]{article}

\usepackage{mathrsfs}
\usepackage{amssymb}
\usepackage{amsmath}
\usepackage{bm}
\usepackage{graphicx}
\usepackage{pst-poly}  
\usepackage{pst-plot}  

\def \la{\lambda}
\def \over{\overline}

\newtheorem{theorem}{Theorem}[section]
\newtheorem{lemma}[theorem]{Lemma}
\newtheorem{coro}[theorem]{Corollary}

\newtheorem{defn}[theorem]{Definition}
\newtheorem{exam}[theorem]{Example}
\newtheorem{rem}[theorem]{Remark}

\newenvironment{pf}[1][Proof]{\noindent\textbf{#1.} }{\hfill\rule{1mm}{2mm}}

\makeatletter \@addtoreset{equation}{section} \makeatother

\begin{document}

\title
{Vulnerability of super edge-connected graphs\thanks{The work was
supported by NNSF of China (No. 61272008).} }

\author
{Zhen-Mu Hong\footnote{E-mail address: zmhong@mail.ustc.edu.cn
(Z.-M. Hong)} \quad  Jun-Ming Xu\footnote{Corresponding author,
E-mail address: xujm@ustc.edu.cn (J.-M. Xu)}
\\ \\
{\small School of Mathematical Sciences}\\
{\small University of Science and Technology of China}\\
{\small  Wentsun Wu Key Laboratory of CAS}  \\
{\small Hefei, Anhui, 230026, China}}

\date{} \maketitle

\begin{center}
\begin{minipage}{140mm}

\begin{center} {\bf Abstract} \end{center}

A subset $F$ of edges in a connected graph $G$ is a $h$-extra
edge-cut if $G-F$ is disconnected and every component has more than
$h$ vertices. The $h$-extra edge-connectivity $\la^{(h)}(G)$ of $G$
is defined as the minimum cardinality over all $h$-extra edge-cuts
of $G$. A graph $G$, if $\la^{(h)}(G)$ exists, is super-$\la^{(h)}$
if every minimum $h$-extra edge-cut of $G$ isolates at least one
connected subgraph of order $h+1$. The persistence $\rho^{(h)}(G)$
of a super-$\la^{(h)}$ graph $G$ is the maximum integer $m$ for
which $G-F$ is still super-$\la^{(h)}$ for any set $F\subseteq E(G)$
with $|F|\leqslant m$. Hong {\it et al.} [Discrete Appl. Math. 160
(2012), 579-587] showed that
$\min\{\la^{(1)}(G)-\delta(G)-1,\delta(G)-1\}\leqslant
\rho^{(0)}(G)\leqslant \delta(G)-1$, where $\delta(G)$ is the
minimum vertex-degree of $G$. This paper shows that
$\min\{\la^{(2)}(G)-\xi(G)-1,\delta(G)-1\}\leqslant \rho^{(1)}(G)\leqslant
\delta(G)-1$, where $\xi(G)$ is the minimum edge-degree of $G$. In
particular, for a $k$-regular super-$\la'$ graph $G$,
$\rho^{(1)}(G)=k-1$ if $\la^{(2)}(G)$ does not exist or $G$ is
super-$\la^{(2)}$ and triangle-free, from which the exact values of
$\rho^{(1)}(G)$ are determined for some well-known networks.

\vskip6pt



\end{minipage}
\end{center}

\vskip0.6cm

\newpage

\section{Introduction}

We follow~\cite{x03} for graph-theoretical terminology and notation
not defined here. Let $G=(V, E)$ be a simple connected graph, where
$V=V(G)$ is the vertex-set of $G$ and $E=E(G)$ is the edge-set of
$G$. It is well known that when the underlying topology of an
interconnection network is modeled by a connected graph $G=(V,E)$,
where $V$ is the set of processors and $E$ is the set of
communication links in the network, the edge-connectivity $\la(G)$
of $G$ is an important measurement for reliability and fault
tolerance of the network. In general, the larger $\la(G)$ is, the
more reliable a network is. Because the connectivity has some
shortcomings, F\`abrega and Fiol~\cite{ff94,ff96} generalized the
concept of the edge-connectivity to the $h$-extra edge-connectivity
for a graph.

\begin{defn}{\rm
Let $h\geqslant 0$ be an integer. A subset $F\subseteq E(G)$ is an
{\it $h$-extra edge-cut} if $G-F$ is disconnected and every
component of $G-F$ has more than $h$ vertices. The {\it $h$-extra
edge-connectivity} of $G$, denoted by $\la^{(h)}(G)$, is defined as
the minimum cardinality of an $h$-extra edge-cut of $G$. }
\end{defn}

Clearly, $\la^{(0)}(G)=\la(G)$ and $\la^{(1)}(G)=\la'(G)$ for any
graph $G$, the latter is called the restricted edge-connectivity
proposed by Esfahanian and Hakimi~\cite{eh88}, who proved that for a
connected graph $G$ of order at least $4$, $\la'(G)$ exists if and
only if $G$ is not a star.

In general, $\la^{(h)}(G)$ does not always exist for $h\geqslant 1$.
For example, let $G^*_{n,h}$ ($n\geqslant h$) be a graph obtained
from $n$ copies of a complete graph $K_h$ of order $h$ by adding a
new vertex $x$ and linking $x$ to every vertex in each of $n$
copies. Clearly, $G^*_{n,1}$ is a star $K_{1,n}$. It is easy to
check that $\la^{(h)}(G^*_{n,h})$ does not exists for $h\geqslant
1$.

A graph $G$ is said a {\it $\la^{(h)}$-graph} or to be {\it
$\la^{(h)}$-connected} if $\la^{(h)}(G)$ exists, and to be {\it not
$\la^{(h)}$-connected} otherwise. For a $\la^{(h)}$-graph $G$, an
$h$-extra edge-cut $F$ is a {\it $\la^{(h)}$-cut} if
$|F|=\la^{(h)}(G)$. It is easy to verify that, for a
$\lambda^{(h)}$-graph $G$,
 \begin{equation}\label{e1.1}
 \lambda^{(0)}(G)\leqslant \lambda^{(1)}(G)\leqslant \lambda^{(2)}(G)\leqslant
 \cdots \leqslant  \lambda^{(h-1)}(G)\leqslant \lambda^{(h)}(G).
 \end{equation}

For two disjoint subsets $X$ and $Y$ in $V(G)$, use $[X,Y]$ to
denote the set of edges between $X$ and $Y$ in $G$. In particular,
$E_G(X)=[X,\overline X]$ and let $d_G(X)=|E_G(X)|$, where $\overline
X=V(G)\setminus X$. For a $\la^{(h)}$-graph $G$, there is certainly
a subset $X\subset V(G)$ with $|X|\geqslant h+1$ such that $E_G(X)$
is a $\la^{(h)}$-cut and, both $G[X]$ and $G[\over{X}]$ are
connected. Such an $X$ is called a {\it $\la^{(h)}$-fragment} of $G$.

For a subset $X\subset V(G)$, use $G[X]$ to denote the subgraph of $G$
induced by $X$. Let
 $$
 \xi_h(G)=\min\{d_G(X):\ X\subset V(G),\
 |X|=h+1\ {\rm and}\ G[X] \ {\rm is\ connected}\}.
 $$
Clearly, $\xi_0(G)=\delta(G)$, the minimum vertex-degree of $G$, and
$\xi_1(G)=\xi(G)$, the minimum edge-degree of $G$ defined as
$\min\{d_G(x)+d_G(y)-2: xy\in E(G)\}$. For a $\lambda^{(h)}$-graph
$G$, Whitney's inequality shows $\lambda^{(0)}(G)\leqslant\xi_0(G)$;
Esfahanian and Hakimi~\cite{eh88} showed
$\lambda^{(1)}(G)\leqslant\xi_1(G)$; Bonsma {\it et
al.}~\cite{buv02}, Meng and Ji~\cite{mj02} showed
$\lambda^{(2)}(G)\leqslant\xi_2(G)$. For $h\geqslant 3$, Bonsma {\it
et al.}~\cite{buv02} found  that the inequality
$\lambda^{(h)}(G)\leqslant\xi_h(G)$ is no longer true in general.
The following theorem shows existence of $\la^{(h)}(G)$ for any
graph $G$ with $\delta(G)\geqslant h$ except for $G^*_{n,h}$.

\begin{theorem}{\rm (Zhang and Yuan~\cite{zy05})}\label{thm1.2}
Let $G$ be a connected graph with order at least $2(\delta+1)$,
where $\delta=\delta(G)$. If $G$ is not isomorphic to
$G^*_{n,\delta}$, then $\la^{(h)}(G)$ exists and
$$
\la^{(h)}(G)\leqslant \xi_h(G)\ \ {\rm for\ any}\ h\ {\rm with}\
0\leqslant h\leqslant \delta.
$$
\end{theorem}

A graph $G$ is said to be {\it $\la^{(h)}$-optimal} if
$\la^{(h)}(G)=\xi_h(G)$. In view of practice in networks, it seems
that the larger $\la^{(h)}(G)$ is, the more reliable the network is.
Thus, investigating $\la^{(h)}$-optimal property of networks has
attracted considerable research interest (see Xu~\cite{x01}). A
stronger concept than $\la^{(h)}$-optimal is super-$\la^{(h)}$.

\begin{defn}{\rm
A $\la^{(h)}$-optimal graph $G$ is {\it super $k$-extra
edge-connected} ({\it super-$\la^{(h)}$} for short), if every
$\la^{(h)}$-cut of $G$ isolates at least one connected subgraph of
order $h+1$. }
\end{defn}

By definition, a super-$\la^{(h)}$ graph is certainly
$\la^{(h)}$-optimal, but the converse is not true. For example, a
cycle of length $n\,(n\geqslant 2h+4)$ is a $\la^{(h)}$-optimal
graph and not super-$\la^{(h)}$. The following necessary and
sufficient condition for a graph to be super-$\la^{(h)}$ is simple
but very useful.

\begin{lemma}\label{lem1.4}
A $\la^{(h)}$-graph $G$ is super-$\la^{(h)}$ if and only if either
$G$ is not $\la^{(h+1)}$-connected or $\la^{(h+1)}(G)>\xi_h(G)$ for
any $h\geqslant 0$.
\end{lemma}

Faults of some communication lines in a large-scale system are
inevitable. However, the presence of faults certainly affects the
super connectedness. The following concept is proposed naturally.

\begin{defn}{\rm
The {\it persistence} of a super-$\la^{(h)}$ graph $G$, denoted by
$\rho^{(h)}(G)$, is the maximum integer $m$ for which $G-F$ is still
super-$\la^{(h)}$ for any subset $F\subseteq E(G)$ with $|F|\leqslant
m$.}
\end{defn}

It is clear that the persistence $\rho^{(h)}(G)$ is a measurement
for vulnerability of super-$\la^{(h)}$ graphs. We can easily obtain
an upper bound on $\rho^{(h)}(G)$ as follows.

\begin{theorem}\label{thm1.6}
$\rho^{(h)}(G)\leqslant \delta(G)-1$ for any super-$\la^{(h)}$ graph
$G$.
\end{theorem}

\begin{pf}
Let $G$ be a super-$\la^{(h)}$ graph and $F$ a set of edges incident
with some vertex of degree $\delta(G)$. Since $G-F$ is disconnected,
$G-F$ is not super-$\la^{(h)}$. By the definition of
$\rho^{(h)}(G)$, we have $\rho^{(h)}(G)\leqslant \delta(G)-1$.
\end{pf}

\vskip6pt

By Theorem~\ref{thm1.6}, we can assume $\delta(G)\geqslant 2$ when
we consider $\rho^{(h)}(G)$ for a super-$\la^{(h)}$ graph $G$. In
this paper, we only focus on the lower bound on $\rho^{(1)}(G)$ for
a super-$\la^{(1)}$ graph $G$. For convenience, we write
$\la^{(0)}$, $\la^{(1)}$, $\la^{(2)}$, $\rho^{(0)}$ and $\rho^{(1)}$
for $\la$, $\la'$, $\la''$, $\rho$ and $\rho'$, respectively.

Very recently, Hong, Meng and Zhang~\cite{hmz12} have showed
$\rho(G)\geqslant \min\{\la'(G)-\delta(G)-1, \delta(G)-1\}$ for any
super-$\la$ and $\la'$-graph $G$. In this paper, we establish
$\rho'(G)\geqslant\min\{\la''(G)-\xi(G)-1,\delta(G)-1\}$,
particularly, for a $k$-regular super-$\la'$ graph $G$,
$\rho'(G)=k-1$ if $G$ is not $\la''$-connected or super-$\la''$ and
triangle-free. As applications, we determine the exact values of
$\rho'$ for some well-known networks.

The left of this paper is organized as follows. In Section 2, we
establish the lower bounds on $\rho'$ for general super-$\la'$
graphs. In Section 3, we focus on regular graphs and give some
sufficient conditions under which $\rho'$ reaches its upper bound or
the difference between upper and lower bounds is at most one. In
Section 4, we determine exact values of $\rho'$ for two well-known
families of networks.

\section{Lower bounds on $\bm{\rho'}$ for general graphs}

In this section, we will establish some lower bounds on $\rho'$ for
a general super-$\la'$ graph. The following lemma is useful for the
proofs of our results.

\begin{lemma}{\rm (Hellwig and Volkmann~\cite{hv05})}\label{lem2.1}
If $G$ is a $\la'$-optimal graph, then $\la(G) = \delta(G)$.
\end{lemma}

\begin{lemma}\label{lem2.2}
Let $G$ be a $\la'$-graph and $F$ be any subset of $E(G)$.

\,{\rm (i)}\ \ If $G$ is $\la'$-optimal and $|F|\leqslant
\delta(G)-1$, then $G-F$ is $\la'$-connected.

{\rm (ii)}\ If $G-F$ is $\la''$-connected, then $G$ is also
$\la''$-connected. Moreover,
\begin{equation}\label{e2.1}
\la''(G-F)\geqslant \la''(G)-|F|.
\end{equation}
\end{lemma}

\begin{pf}
Let $G$ be a $\la'$-graph of order $n$ and $F$ be any subset of $E(G)$.
Clearly, $n\geqslant 4$.

(i) Assume that $G$ is $\la'$-optimal and $|F|\leqslant \delta(G)-1$.
It is trivial for $\delta(G)=1$. Assume $\delta(G)\geqslant 2$ below.
Since $G$ is $\la'$-optimal, $\la(G)=\delta(G)$ by Lemma
\ref{lem2.1}. By $|F|\leqslant \delta(G)-1$, $G-F$ is connected. If
$G-F$ is a star $K_{1,n-1}$, then $G$ has a vertex $x$ with degree
$n-1$. Let $H=G-x$. Then $F=E(H)$ and $\delta(H)\geqslant
\delta(G)-1$. Thus,
 $$
 \delta(G)-1\geqslant |F|=|E(H)|\geqslant \frac12(n-1)(\delta(G)-1),
 $$
which implies $n\leqslant 3$, a contradiction. Thus, $G-F$ is not a
star $K_{1,n-1}$, and so is $\la'$-connected.

(ii) Assume that $G-F$ is $\la''$-connected, and let $X$ be a
$\la''$-fragment of $G-F$. Clearly, $E_G(X)$ is a $2$-extra edge-cut
of $G$, and so $G$ is $\la''$-connected and $d_G(X)\geqslant
\la''(G)$. Thus, $\la''(G-F)=d_{G-F}(X) \geqslant
d_G(X)-|F|\geqslant \la''(G)-|F|.$
\end{pf}

\vskip6pt

By Lemma~\ref{lem2.2}, we obtain the following result immediately.

\begin{theorem}\label{thm2.3}
Let $G$ be a super-$\la'$ graph. If $G$ is not $\la''$-connected,
then $\rho'(G)=\delta(G)-1$.
\end{theorem}

\begin{pf}
Since $G$ is super-$\la'$, $G$ is $\la'$-optimal. Let $F$
be any subset of $E(G)$ with $|F|\leqslant \delta(G)-1$. By
Lemma~\ref{lem2.2} (i), $G-F$ is $\la'$-connected. If $G$ is not
$\la''$-connected, then $G-F$ is also not $\la''$-connected by
Lemma~\ref{lem2.2} (ii). By Lemma~\ref{lem1.4} $G-F$ is
super-$\la'$, which implies $\rho'(G)\geqslant \delta(G)-1$.
Combining this with Theorem~\ref{thm1.6}, we obtain the conclusion.
\end{pf}

\vskip6pt

By Theorem~\ref{thm2.3}, we only need to consider $\rho'(G)$ for a
$\la''$-connected super-$\la'$ graph $G$. A graph $G$ is said to be
{\it edge-regular} if $d_G(\{x,y\})=\xi(G)$ for every $xy\in E(G)$,
where $d_G(\{x,y\})$ is called the edge-degree of the edge $xy$ in
$G$. Denote by $\eta(G)$ the number of edges with edge-degree
$\xi(G)$ in $G$. For simplicity, we write $\la''=\la''(G)$,
$\la'=\la'(G)$, $\rho'=\rho'(G)$, $\xi=\xi(G)$ and
$\delta=\delta(G)$ when just one graph $G$ is under discussion.

\begin{theorem}\label{thm2.4}
Let $G$ be a $\la''$-connected super-$\la'$ graph. Then

\, {\rm (i)}\ \ $\rho'(G)\geqslant \min\{\la''-\xi-1,\delta-1\}$ if
$\eta(G)\geqslant \delta$, or

{\rm (ii)}\ $\rho'(G)\geqslant \min\{\la''-\xi,\delta-1\}$ if $G$ is
edge-regular.

\end{theorem}

\begin{pf}
Since $G$ is $\la''$-connected and super-$\la'$, $\la''>\xi$ by
Lemma~\ref{lem1.4}. If $\delta=1$, then $\rho'=0$. Assume
$\delta\geqslant 2$ below. Let
 \begin{equation}\label{e2.2}
 \begin{array}{ll}
 & m_1 =\min\{\la''-\xi-1,\delta-1\}, \\
 & m_2 =\min\{\la''-\xi,\delta-1\} \ \ {\rm and}\ \ m=m_1\ \ {\rm or}\ \ m_2.
 \end{array}
 \end{equation}
Since $\la''>\xi$ and $\delta\geqslant 2$, $0\leqslant m_1\leqslant \la''-\xi-1$,
$1\leqslant m_2\leqslant \la''-\xi$ and $m\leqslant\delta-1$.

Let $F$ be any subset of $E(G)$ with $|F|= m$ and let $G'=G-F$. Since
$G$ is $\la'$-optimal and $|F|\leqslant \delta-1$, $G'$ is $\la'$-connected by
Lemma~\ref{lem2.2} (i). To show that $\rho'\geqslant m$, we only need to
prove that $G'$ is super-$\la'$. If $G'$ is not $\la''$-connected,
then $G'$ is super-$\la'$ by Lemma~\ref{lem1.4}. Assume now that
$G'$ is $\la''$-connected. It follows from (\ref{e2.1}) and (\ref{e2.2}) that
 \begin{equation}\label{e2.3}
 \la''(G')\geqslant\la''(G)-|F|= \la''-m \geqslant
  \left\{\begin{array}{ll}
 \xi+1 & {\rm if}\ m=m_1;\\
 \xi & {\rm if}\ m=m_2.
  \end{array}\right.
 \end{equation}

Since $|F|\leqslant \delta-1$, if $\eta(G)\geqslant \delta$, $G'$
has at least one edge with edge-degree $\xi(G)$, which implies
$\xi(G')\leqslant \xi(G)$. Moreover, if $G$ is edge-regular, then
$\eta(G)\geqslant \delta$ and every edge of $G$ is incident with
some edge with edge-degree $\xi$, which implies $\xi(G')< \xi(G)$ if
$|F|\geqslant 1$. It follows that
 \begin{equation}\label{e2.4}
 \xi(G')\leqslant \left\{\begin{array}{ll}
 \xi(G) & {\rm if}\ \eta(G)\geqslant \delta;\\
 \xi(G)-1 & {\rm if}\ G\ \text{is\ edge-regular\ and $|F|=m\geqslant 1$}.
 \end{array}\right.
\end{equation}

Combining (\ref{e2.3}) with (\ref{e2.4}), if $m=m_1$
and $\eta(G)\geqslant \delta$ or $m=m_2\geqslant 1$
and $G$ is edge-regular, we have $\la''(G')>\xi(G')$. By Lemma
\ref{lem1.4}, $G'$ is super-$\la'$, and so the conclusions (i) and
(ii) hold.

The theorem follows.
\end{pf}

\begin{figure}[h]
\begin{center}
\begin{pspicture}(-2.5,0)(2.5,3.0)

\cnode(-2.5,2.1){3pt}{a} \cnode(-2.5,0.9){3pt}{b}\ncline{a}{b}
\cnode(-1.5,3){3pt}{c} \cnode(-1.5,0){3pt}{d}\ncline{c}{d}
\cnode(-.3,3){3pt}{e} \cnode(-.3,0){3pt}{f}\ncline{e}{f}
\cnode(0.7,2.1){3pt}{g} \cnode(0.7,0.9){3pt}{h}\ncline{g}{h}
\cnode(2,2.1){3pt}{i} \cnode(2,0.9){3pt}{j}\ncline{i}{j}
\rput(2.2,1.5){$e$}
\ncline{a}{c}\ncline{a}{d}\ncline{a}{g}
\ncline{b}{c}\ncline{b}{d}\ncline{b}{h}
\ncline{c}{e}\ncline{d}{f}
\ncline{g}{e}\ncline{g}{f}\ncline{g}{i}
\ncline{h}{e}\ncline{h}{f}\ncline{h}{j}

\end{pspicture}
\caption{{\small The graph $G$ in Remark
\ref{rem2.5}.}\label{fig1}}
\end{center}
\end{figure}
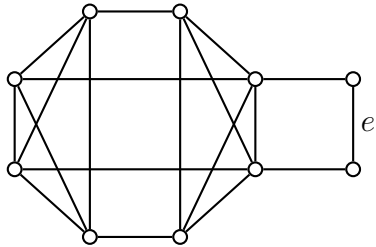

\begin{rem}\label{rem2.5}
{\rm The condition ``\,$\eta(G)\geqslant \delta$\," in Theorem
\ref{thm2.4} is necessary. For example, consider the graph $G$ shown
in Figure \ref{fig1}, $\eta(G)=1<2=\delta$. Since
$\xi=\la=\delta=2<4=\la''$, $G$ is super-$\la'$ by
Lemma~\ref{lem1.4}. We should have that $\rho'(G)\geqslant
\la''-\xi-1=\delta-1=1$ by Theorem~\ref{thm2.4}, which shows the
removal of any edge from $G$ results in a super-$\la'$ graph.
However, $\la''(G-e)=\xi(G-e)=4$, and so $G-e$ is not super-$\la'$
by Lemma~\ref{lem1.4}, which implies $\rho'(G)=0$.}
\end{rem}

The {\it Cartesian product} of graphs $G_1$ and $G_2$ is the graph
$G_1\times G_2$ with vertex-set $V(G_1)\times V(G_2)$, two vertices
$x_1x_2$ and $y_1y_2$, where $x_1,y_1\in V(G_1)$ and $x_2,y_2\in
V(G_2)$, being adjacent in $G_1\times G_2$ if and only if either
$x_1=y_1$ and $x_2y_2\in E(G_2)$, or $x_2=y_2$ and $x_1y_1\in
E(G_1)$. The study on $\la'$ for Cartesian products can be found in
\cite{lcm09,lcx07,ou11}.

\begin{figure}[ht]
\begin{center}
\hspace*{30pt}
\begin{pspicture}(-9.6,-1.5)(3.5,4.8)

\cnode(-9.5,4.3){3pt}{a1} \cnode(-7,4.3){3pt}{b1} \ncline{a1}{b1}
\cnode(-9.5,3.0){3pt}{a2} \cnode(-7,3.0){3pt}{b2} \ncline{a2}{b2}
\cnode(-9.5,1.5){3pt}{a3} \cnode(-7,1.5){3pt}{b3} \ncline{a3}{b3}
\cnode(-9.5,0.7){3pt}{a4} \cnode(-7,0.7){3pt}{b4} \ncline{a4}{b4}
\cnode(-9.5,-.3){3pt}{a5} \cnode(-7,-.3){3pt}{b5} \ncline{a5}{b5}
\cnode(-6.1,3.6){3pt}{x1} \cnode(-6.1,2.35){3pt}{y1}\ncline{x1}{y1}
\cnode(-5.2,3.6){3pt}{x2} \cnode(-5.2,2.35){3pt}{y2}\ncline{x2}{y2}
\cnode(-3.5,3.6){3pt}{x3} \cnode(-3.5,2.35){3pt}{y3}\ncline{x3}{y3}
\ncline{x1}{b1}\ncline{x1}{b2}\ncline{x1}{b3}
\ncline{y1}{b1}\ncline{y1}{b2}\ncline{y1}{b3}
\ncline{x2}{b1}\ncline{x2}{b2}\ncline{x2}{b3}
\ncline{y2}{b1}\ncline{y2}{b2}\ncline{y2}{b3}
\ncline{x3}{b1}\ncline{x3}{b2}\ncline{x3}{b3}
\ncline{y3}{b1}\ncline{y3}{b2}\ncline{y3}{b3}
\rput(-5.84,3.52){\scriptsize $x_{_1}$}\rput(-5.84,2.36){\scriptsize $y_{_1}$}
\rput(-4.9,3.6){\scriptsize $x_{_2}$}\rput(-4.9,2.35){\scriptsize $y_{_2}$}
\rput(-3.2,3.6){\scriptsize $x_{_t}$}\rput(-3.2,2.35){\scriptsize $y_{_t}$}

\rput(-9.5,2.4){$\vdots$}\rput(-7,2.4){$\vdots$}
\rput(-9.5,0.3){$\vdots$}\rput(-7,0.3){$\vdots$}
\rput(-4.35,3.6){$\cdots$}\rput(-4.35,2.35){$\cdots$}
\rput(-7,1.05){\rnode{C1}{}} \rput(-7,4.5){\rnode{C2}{}}
\ncbox[nodesep=3pt,boxsize=.4]{C1}{C2} \rput(-7,1.18){\small $Z$}
\rput(-9.5,-0.9){\rnode{A1}{}} \rput(-9.5,4.7){\rnode{A2}{}}
\ncbox[nodesep=3pt,boxsize=.6]{A1}{A2} \rput(-9.5,-0.7){\small $X$}
\rput(-7,-0.9){\rnode{B1}{}} \rput(-7,4.7){\rnode{B2}{}}
\ncbox[nodesep=3pt,boxsize=.6]{B1}{B2} \rput(-7,-0.7){\small $Y$}
\rput(-8.25,-1.4){$G$}

\cnode(-1,4.5){3pt}{p1}\cnode(0.8,4.5){3pt}{p2}\cnode(2.6,4.5){3pt}{p3}
\ncline{p1}{p2}\ncline{p2}{p3} \pscurve[curvature=2 0.1 0](-0.94,4.55)(0.8,4.8)(2.54,4.55)
\cnode(-1,3.7){3pt}{p4}\cnode(0.8,3.7){3pt}{p5}\cnode(2.6,3.7){3pt}{p6}
\ncline{p4}{p5}\ncline{p5}{p6} \pscurve[curvature=2 0.1 0](-0.94,3.75)(0.8,4)(2.54,3.75)
\cnode(-1,2.9){3pt}{p7}\cnode(0.8,2.9){3pt}{p8}\cnode(2.6,2.9){3pt}{p9}
\ncline{p7}{p8}\ncline{p8}{p9} \pscurve[curvature=2 0.1 0](-0.94,2.95)(0.8,3.2)(2.54,2.95)
\ncline{p1}{p4}\ncline{p4}{p7} \pscurve[curvature=2 0.1 0](-1.06,4.45)(-1.27,3.7)(-1.06,2.95)
\ncline{p2}{p5}\ncline{p5}{p8} \pscurve[curvature=2 0.1 0](0.74,4.45)(0.53,3.7)(0.74,2.95)
\ncline{p3}{p6}\ncline{p6}{p9} \pscurve[curvature=2 0.1 0](2.54,4.45)(2.33,3.7)(2.54,2.95)

\cnode(-1,1.6){3pt}{q1}\cnode(0.8,1.6){3pt}{q2}\cnode(2.6,1.6){3pt}{q3}
\ncline{q1}{q2}\ncline{q2}{q3} \pscurve[curvature=2 0.1 0](-0.94,1.55)(0.8,1.3)(2.54,1.55)
\cnode(-1,0.8){3pt}{q4}\cnode(0.8,0.8){3pt}{q5}\cnode(2.6,0.8){3pt}{q6}
\ncline{q4}{q5}\ncline{q5}{q6} \pscurve[curvature=2 0.1 0](-0.94,0.75)(0.8,0.5)(2.54,0.75)
\cnode(-1,0){3pt}{q7}\cnode(0.8,0){3pt}{q8}\cnode(2.6,0){3pt}{q9}
\ncline{q7}{q8}\ncline{q8}{q9} \pscurve[curvature=2 0.1 0](-0.94,-0.05)(0.8,-0.3)(2.54,-0.05)
\ncline{q1}{q4}\ncline{q4}{q7} \pscurve[curvature=2 0.1 0](-1.06,1.55)(-1.27,0.8)(-1.06,0.05)
\ncline{q2}{q5}\ncline{q5}{q8} \pscurve[curvature=2 0.1 0](0.74,1.55)(0.53,0.8)(0.74,0.05)
\ncline{q3}{q6}\ncline{q6}{q9} \pscurve[curvature=2 0.1 0](2.54,1.55)(2.33,0.8)(2.54,0.05)

\ncline[linewidth=0.6mm]{p7}{q1}\ncline{p8}{q2}\ncline{p9}{q3}
\pscurve[curvature=2 0.1 0](-0.94,3.65)(-0.7,2.25)(-0.94,0.85)   
\pscurve[curvature=2 0.1 0](0.86,3.65)(1.1,2.25)(0.86,0.85)      
\pscurve[curvature=2 0.1 0](2.66,3.65)(2.9,2.25)(2.66,0.85)      
\pscurve[curvature=2 0.1 0](-0.94,4.45)(-0.4,2.25)(-0.94,0.05)   
\pscurve[curvature=2 0.1 0](0.86,4.45)(1.4,2.25)(0.86,0.05)      
\pscurve[linewidth=0.6mm,curvature=2 0.1 0](2.66,4.45)(3.2,2.25)(2.66,0.05)      
\rput(-1.2,2.25){\small $e_{_1}$}\rput(3.45,2.25){\small $e_{_2}$}
\rput(0.8,-1.2){$H=K_2\times K_3\times K_3$}
\end{pspicture}
\caption{\small Two graphs $G$ and $H$ in Remark
\ref{rem2.7}.\label{fig2}}
\end{center}
\end{figure}
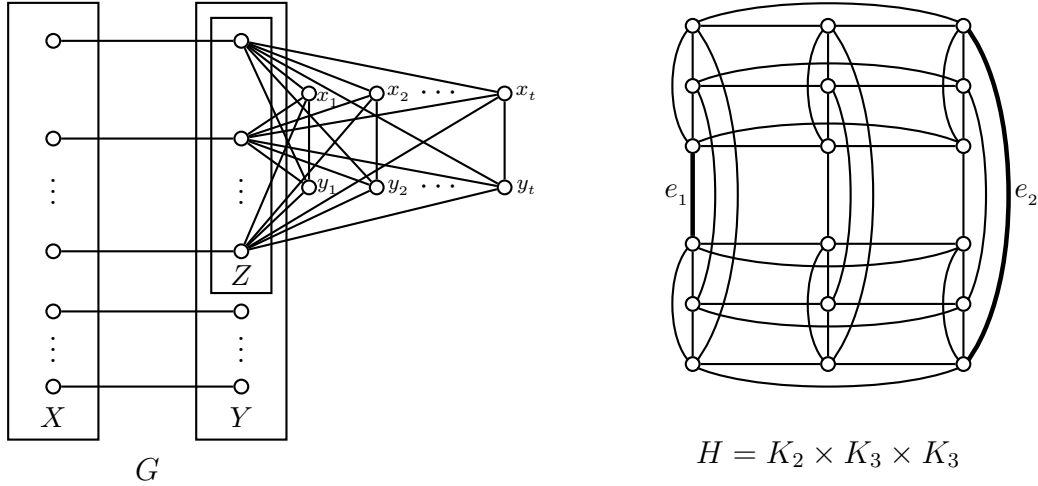

\begin{rem}\label{rem2.7}
{\rm The graphs $G$ and $H$ shown in Figure \ref{fig2} can show that
the lower bounds on $\rho'$ given in Theorem~\ref{thm2.4} are sharp.

In $G$, $X$ and $Y$ are two disjoint subsets of $3t-2$ vertices, and
$Z$ is a subset of $Y$ with $t-1$ vertices, where $t\geqslant 2$.
There is a perfect matching between $X$ and $Y$ and the subgraphs
induced by $X,Y$ and $Z\cup\{x_i,y_i\}$ are all complete graphs, for
each $i=1,2,\ldots,t$. It is easy to check that
$\eta(G)=\delta(G)=\la(G)=t$, $\la'(G)=\xi(G)=2t-2$, $\la''(G)=3t-2$
and $G$ is super-$\la'$. By Theorem \ref{thm2.4}, $\rho'(G)\geqslant
\la''(G)-\xi(G)-1=t-1=\delta(G)-1$. Combining this fact with Theorem
\ref{thm1.6}, we have $\rho'(G)=\delta(G)-1$. This example shows
that the lower bound on $\rho'$ given in Theorem~\ref{thm2.4} (i) is
sharp.

For the 5-regular graph $H=K_2\times K_3\times K_3$, $\la''(H)=9$
and $\xi(H)=8$, and so $H$ is super-$\la'$ by Lemma \ref{lem1.4}. On
the one hand, by Theorem \ref{thm2.4}, $\rho'(H)\geqslant
\la''(H)-\xi(H)=1$. On the other hand, for $F=\{e_1,e_2\}$,
$\la''(H-F)=7=\xi(H-F)$, and so $H-F$ is not super-$\la'$ by
Lemma~\ref{lem1.4}, which yields $\rho'(H)\leqslant 1$. Hence,
$\rho'(H)=\la''(H)-\xi(H)=1$. This example shows that the lower
bound on $\rho'$ given in Theorem~\ref{thm2.4} (ii) is sharp.}
\end{rem}

\section{Bounds on $\bm{\rho'}$ for regular graphs}

The {\it girth}  of a graph $G$, denoted by $g(G)$, is the length of
a shortest cycle in $G$. A graph is said to be {\it $C_n$-free} if
it contains no cycles of length $n$. In general, $C_3$-free is said
{\it triangle-free}. To guarantee that $G$ is edge-regular, which is
convenient for us to use Theorem~\ref{thm2.4}, we consider regular
graphs in this section.

Clearly, any $k$-regular graph contains cycles if $k\geqslant 2$. It
is easy to check that $C_4$ and $C_5$ are only two $2$-regular
super-$\la'$ graphs. Obviously, $\rho'(C_4)=\rho'(C_5)=1$. In the
following discussion, we always assume $k\geqslant 3$ when we
mention $k$-regular connected graphs. We first consider $3$-regular
graphs, such graphs have even order.

\begin{lemma}\label{lem3.1}
Let $G$ be a $3$-regular super-$\la'$ graph of order $2n$. If
$n\geqslant 4$, then the grith $g(G)>4$ and $n\ne 4$.
\end{lemma}

\begin{pf}
Since $G$ is a $3$-regular super-$\la'$ graph of order at least $8$,
$\la'(G)=\xi(G)=4$, and so $\la(G)=3$ by Lemma~\ref{lem2.1}.
Moreover, every $\la'$-cut of $G$ isolates one edge. If $G$ contains
a $C_3$, then let $X=V(C_3)$. If $G-E_G(X)$ isolates a vertex, then
$G\cong K_4$, a contradiction with $n\geqslant 4$. Thus, $E_G(X)$ is
a 1-extra edge-cut and $\la'(G)\leqslant d_G(X)=3<4=\la'(G)$, a
contradiction. If $G$ contains a $C_4$, let $Y=V(C_4)$, then
$G-E_G(Y)$ does not isolate a vertex since $G$ contains no
triangles, and so $E_G(Y)$ is also a 1-extra edge-cut and
$4=\xi(G)=\la'(G)\leqslant d_G(Y)=4$, which implies that $E_G(Y)$ is
$\la'$-cut of $G$ and does not isolate one edge since $n\geqslant
4$, which means that $G$ is not super-$\la'$, a contradiction. Thus,
the girth $g(G)>4$. Moreover, since any $3$-regular graph with girth
greater than 4 has at least 10 vertices, we have $n\geqslant 5$.
\end{pf}

\begin{theorem}\label{thm3.2}
Let $G$ be a $3$-regular super-$\la'$ graph of order $2n$. If $n=2$
or $3$, then $\rho'(G)=2$. If $n\geqslant 5$, then $\rho'(G)=1$.
\end{theorem}

\begin{pf}
The complete graph $K_4$ and the complete bipartite graph $K_{3,3}$
are the unique $3$-regular super-$\la'$ graphs of order $4$ and $6$, respectively.
It is easy to check that $\rho'(K_4)=\rho'(K_{3,3})=2$.

Next, assume $n\geqslant 4$. Then $g(G)>4$ and $n\geqslant 5$ by
Lemma~\ref{lem3.1}. Since $G$ is 3-regular super-$\la'$,
$\la'(G)=\xi(G)=4$ and every $\la'$-cut isolates at least one edge.
Since $g(G)\geqslant 5$, $G$ is not isomorphic to $G^*_{n,2}$. By
Theorem~\ref{thm1.2}, $G$ is $\la''$-connected. By
Lemma~\ref{lem1.4} and Theorem \ref{thm2.4} (ii), $\rho'(G)\geqslant
\la''-\xi\geqslant 1$. To prove $\rho'(G)\leqslant 1$, we only need
to show that there exists a subset $F\subset E(G)$ with $|F|=2$ such
that $G-F$ is not super-$\la'$.

Let $P=(u,v,w)$ be a path of length two in $G$. Since $g(G)>4$, $u$
and $w$ have only common neighbor $v$. Let $\{u_1,u_2,v\}$ and
$\{w_1,w_2,v\}$ are the sets of neighbors of $u$ and $w$,
respectively. Then either $u_1w_1\not\in E(G)$ or $u_1w_2\not\in
E(G)$ since $g(G)>4$. Assume $u_1w_1\not\in E(G)$ and let
$F=\{uu_1,ww_1\}$. Then $\xi(G-F)=3$. Set $X=V(P)$. Then
$d_{G-F}(X)=3=\xi(G-F)$. Moreover, it is easy to see that
$G[\overline X]$ is connected. Thus, $X$ is a $2$-extra edge-cut of
$G-F$, and so $\la''(G-F)\leqslant d_{G-F}(X)=\xi(G-F)$. By Lemma
\ref{lem1.4}, $G-F$ is not super-$\la'$, which yields
$\rho'(G)\leqslant 1$. Hence, $\rho'(G)=1$, and so the theorem
follows.
\end{pf}

\bigskip
The well-known Peterson graph $G$ is a $3$-regular super-$\la'$
graph with girth $g(G)=5$. By Theorem \ref{thm3.2}, $\rho'(G)=1$.

In general, it is quite difficult to determine the exact value of
$\rho'(G)$ of a $k$-regular super-$\la'$ graph $G$ for $k\geqslant
4$. By Theorem~\ref{thm2.3} for a $k$-regular super-$\la'$ graph
$G$, if $G$ is not $\la''$-connected, then $\rho'(G)=k-1$. Thus, we
only need to consider $k$-regular $\la''$-graphs. For such a graph
$G$, we can establish some bounds on $\rho'(G)$ in terms of $k$.

\begin{lemma}\label{lem3.3}
Let $G$ be $k$-regular $\la''$-optimal graph and $k\geqslant 4$.
Then $G$ is super-$\la'$ if and only if $g(G)\geqslant 4$ or $k\geqslant 5$.
\end{lemma}

\begin{pf}
Let $G$ be a $k$-regular $\la''$-optimal graph and $k\geqslant 4$. Then
$\la''(G)=\xi_2(G)=3k-4>2k-2=\xi(G)$ if and only if $g(G)\geqslant 4$,
and $\la''(G)=\xi_2(G)\geqslant 3k-6>2k-2=\xi(G)$ if and only if $k\geqslant 5$.
Either of two cases shows that $G$ is super-$\la'$ by
Lemma~\ref{lem1.4}.
\end{pf}

\begin{theorem}\label{thm3.4}
Let $G$ be a $k$-regular $\la''$-optimal graph and $k\geqslant 4$.
If $g(G)\geqslant 4$, then
 $$
 k-2\leqslant \rho'(G)\leqslant k-1.
 $$
\end{theorem}

\begin{pf}
Since $G$ is $\la''$-optimal, $G$ is super-$\la'$ by
Lemma~\ref{lem3.3}. Since $G$ is $k$-regular and $g(G)\geqslant 4$,
$\la''=3k-4$ and $\xi=2k-2$. By Theorem~\ref{thm2.4} (ii),
$\rho'(G)\geqslant \la''-\xi=k-2$. By Theorem~\ref{thm1.6},
$\rho'(G)\leqslant k-1$.
\end{pf}

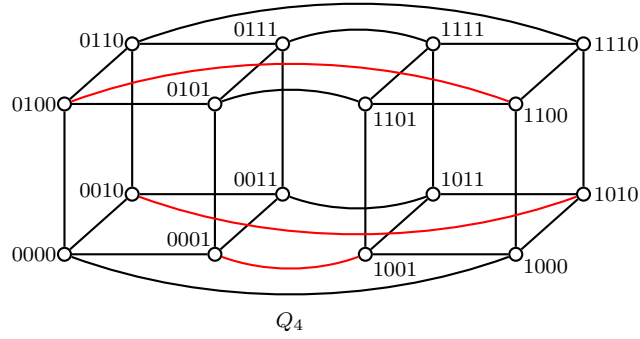
\begin{figure}[h]
\begin{pspicture}(-3,0)(8,4.5)
\cnode(1,1){.1}{0000}\rput(.6,1){\scriptsize0000}
\cnode(1,3){.1}{0100}\rput(.6,3){\scriptsize0100}
\cnode(3,1){.1}{0001}\rput(2.66,1.2){\scriptsize0001}
\cnode(3,3){.1}{0101}\rput(2.66,3.2){\scriptsize0101}
\cnode(1.9,1.8){.1}{0010}\rput(1.5,1.85){\scriptsize0010}
\cnode(1.9,3.8){.1}{0110}\rput(1.5,3.85){\scriptsize0110}
\cnode(3.9,1.8){.1}{0011}\rput(3.55,2){\scriptsize0011}
\cnode(3.9,3.8){.1}{0111}\rput(3.55,4){\scriptsize0111}
\cnode(5,1){.1}{1001}\rput(5.4,.8){\scriptsize1001}
\cnode(5,3){.1}{1101}\rput(5.4,2.8){\scriptsize1101}
\cnode(7,1){.1}{1000}\rput(7.4,.85){\scriptsize1000}
\cnode(7,3){.1}{1100}\rput(7.4,2.85){\scriptsize1100}
\cnode(5.9,1.8){.1}{1011}\rput(6.3,2){\scriptsize1011}
\cnode(5.9,3.8){.1}{1111}\rput(6.3,4){\scriptsize1111}
\cnode(7.9,1.8){.1}{1010}\rput(8.34,1.8){\scriptsize1010}
\cnode(7.9,3.8){.1}{1110}\rput(8.34,3.8){\scriptsize1110}
\ncline{0000}{0001}\ncline{0001}{0101}\ncline{0101}{0100}\ncline{0100}{0000}
\ncline{0010}{0011}\ncline{0011}{0111}\ncline{0111}{0110}\ncline{0110}{0010}
\ncline{0000}{0010}\ncline{0001}{0011}\ncline{0101}{0111}\ncline{0100}{0110}
\ncline{1001}{1000}\ncline{1000}{1010}\ncline{1010}{1011}\ncline{1011}{1001}
\ncline{1101}{1100}\ncline{1100}{1110}\ncline{1110}{1111}\ncline{1111}{1101}
\ncline{1001}{1101}\ncline{1000}{1100}\ncline{1010}{1110}\ncline{1011}{1111}
\nccurve[angleA=-20,angleB=-160]{0000}{1000}
\nccurve[angleA=-20,angleB=-160,linecolor=red]{0010}{1010}
\nccurve[angleA=-20,angleB=-160,linecolor=red]{0001}{1001}
\nccurve[angleA=-20,angleB=-160]{0011}{1011}
\nccurve[angleA=20,angleB=160,linecolor=red]{0100}{1100}
\nccurve[angleA=20,angleB=160]{0110}{1110}
\nccurve[angleA=20,angleB=160]{0101}{1101}
\nccurve[angleA=20,angleB=160]{0111}{1111}
\rput(4,.1){\scriptsize$Q_4$}
\end{pspicture}
\caption{\label{f3} \footnotesize {The hypercube $Q_4$}}
\end{figure}

\begin{rem}\label{rem3.5}\ {\rm
The lower bound on $\rho'$ given in Theorem \ref{thm3.4} is sharp.
For example, the $4$-dimensional cube $Q_4$ (see Figure \ref{f3}) is
a $4$-regular graph with girth $g=4$ and $\la''(Q_4)=\xi_2(Q_4)=8$.
On the one hand, $\rho'(Q_4)\geqslant 2$ by Theorem~\ref{thm3.4}. On
the other hand, let $X$ be the subset of vertices of $Q_4$ whose
first coordinates are 0 and
$F=\{(0001,1001),(0010,1010),(0100,1100)\}$ (shown by red edges in
Figure \ref{f3}). Since $\la''(Q_4-F)\leqslant
d_{Q_4-F}(X)=5=\xi(Q_4-F)$, $Q_4-F$ is not super-$\la'$ by Lemma
\ref{lem1.4}, which implies $\rho'(Q_4)\leqslant 2$. Hence,
$\rho'(Q_4)=2$. }
\end{rem}

\vskip6pt

For a $k$-regular $\la''$-optimal graph with $g(G)=3$, we can establish
an upper bound on $\rho'$ under some conditions. To prove our
result, we need the following lemma.

\begin{lemma}{\rm (Hong {\it et al.}~\cite{hmz12})}\label{lem3.6}
Let $G$ be an $m$-connected graph. Then for any subset $X\subset
V(G)$ with $|X|\geqslant m$ and $|\over{X}| \geqslant m$, there are
at least $m$ independent edges in $E_G(X)$.
\end{lemma}

\begin{theorem}\label{thm3.7}
Let $G$ be a $k$-regular $\la''$-optimal graph with $g(G)=3$ and
$k\geqslant 5$. If $G$ is $(k-2)$-connected and not super-$\la''$,
then
  \begin{equation}\label{e3.1}
  k-4\leqslant \rho'(G)\leqslant k-3,
 \end{equation}
and the bounds are best possible.
\end{theorem}

\begin{pf}
Since $G$ is $k$-regular $\la''$-optimal, $G$ is super-$\la'$ by
Lemma~\ref{lem3.3}, $\la''=3k-6$ and $\xi=2k-2$. By
Theorem~\ref{thm2.4} (ii), $\rho'(G)\geqslant \la''-\xi=k-4$. Thus,
we only need to prove $\rho'(G)\leqslant k-3$.

Since $G$ is not super-$\la''$, there exists a $\la''$-fragment $X$
of $G$ such that $|\overline X|\geqslant |X|\geqslant 4$. Let
$|X|=t$. If $t<k-2$, then $k>6$. For any $x\in X$, since
$d_{G[X]}(x)\leqslant t-1$,
$|[\{x\},\over{X}]|=d_G(x)-d_{G[X]}(x)\geqslant k-t+1$, and so
  \begin{equation}\label{e3.2}
 3k-6=\la''=d_G(X)=\sum_{x\in X}|[\{x\},\over{X}]|\geqslant
t(k-t+1).
 \end{equation}
Since the function $f(t)=t(k-t+1)$ is convex in the integer interval
$[3, k-2]$ and reaches the minimum value at two end-points of the
interval. It follows that
 \begin{equation}\label{e3.3}
 f(t)>f(k-2)=f(3)=3k-6\ \ {\rm for}\ k>6.
 \end{equation}
Comparing (\ref{e3.3}) with (\ref{e3.2}), we obtain a contradiction.
Thus, $t\geqslant k-2$. By Lemma~\ref{lem3.6}, there exists a subset
$F\subseteq E_G(X)$ consisting of $k-2$ independent edges. If $G-F$
is not $\la''$-connected, then $G-F$ is not super-$\la'$ by
Lemma~\ref{lem1.4}. Assume that $G-F$ is $\la''$-connected. Then
$E_{G-F}(X)$ is a $2$-extra edge-cut of $G-F$. Since
 $$
 \la''(G-F)\leqslant d_{G-F}(X)=d_G(X)-|F|=2k-4=\xi-2\leqslant \xi(G-F),
 $$
$G-F$ is not super-$\la'$ by Lemma \ref{lem1.4}. Hence
$\rho'(G)\leqslant k-3$.

To show these bounds are best possible, we consider the graph
$H=K_2\times K_3\times K_3$ and  $G=K_4\times K_4$. For the graph
$H$, it is $5$-regular $\la''$-optimal, and $\rho'(H)=1$ (see Remark
\ref{rem2.7}), which shows that the lower bound given in
(\ref{e3.1}) is sharp when $k=5$. For the graph $G$, it is
$6$-regular $\la''$-optimal but not super-$\la''$. For any subset
$F\subset E(G)$ with $|F|=3$, $G-F$ is certainly $\la''$-connected
and $\la''(G-F)\geqslant \la''-|F|=12-3=9>8\geqslant \xi(G-F)$. By
Lemma~\ref{lem1.4}, $G-F$ is super-$\la'$, which yields
$\rho'(G)\geqslant 3$. Hence, $\rho'(G)=3$, which shows that the
upper bound given in (\ref{e3.1}) is sharp.

The theorem follows.
\end{pf}

\vskip6pt

For a $k$-regular super-$\la''$ graph, the lower bound on $\rho'$
can be improved a little, which is stated as the following theorem.

\begin{theorem}\label{thm3.8}
Let $G$ be a $k$-regular super-$\la''$ graph and $k\geqslant 4$.
Then

 \, {\rm (i)}\ \ $\rho'(G)=k-1$ if $k\geqslant 4$ when $g(G)\geqslant 4$;

 {\rm (ii)}\, $\rho'(G)\geqslant k-3$ if $k\geqslant 6$ and $\rho'(G)=2$ if $k=5$ when $g(G)=3$.
 \end{theorem}

\begin{pf}
Since $G$ is super-$\la''$, $G$ is $\la''$-optimal and
$\la''\leqslant 3k-4$. If $g(G)\geqslant 4$ or $k\geqslant 5$, then
$G$ is super-$\la'$ by Lemma~\ref{lem3.3}. Let $F$ be any subset of
$E(G)$ with $|F|=\la''-\xi+1$ and $G'=G-F$. Since
$|F|=\la''-\xi+1\leqslant k-1$, $G'$ is $\la'$-connected by
Lemma~\ref{lem2.2} (i). We first prove that
 \begin{equation}\label{e3.5}
 \rho'(G)\geqslant \la''-\xi+1\ \text{~if $g(G)\geqslant 4$ or $k\geqslant 5$}.
 \end{equation}
To the end, we need to prove that $G'$ is super-$\la'$. By Lemma
\ref{lem1.4}, we only need to prove that
\begin{equation}\label{e3.6}
\la''(G')>\xi(G')\ \ \text{if $G'$\ is\ $\la''$-connected}.
\end{equation}

Let $X$ be any $\la''$-fragment of $G'$. Since $d_{G'}(\overline
X)=d_{G'}(X)=\la''(G')$, we can assume $|X|\leqslant |\overline X|$.
Since $X$ is a $2$-extra edge-cut of $G$, $d_G(X)\geqslant
\la''(G)=\la''$, and so
\begin{equation}\label{e3.7}
\la''(G')=d_{G'}(X)\geqslant d_G(X)-|F|\geqslant \la''-|F|=\xi-1.
\end{equation}
On the other hand, since $G$ is edge-regular, we have
\begin{equation}\label{e3.8}
\xi(G')\leqslant \xi(G)-1=\xi-1.
\end{equation}
Combing (\ref{e3.7}) with (\ref{e3.8}), in order to prove
(\ref{e3.6}), we only need to show that at least one of the
inequalities (\ref{e3.7}) and (\ref{e3.8}) is strict.

If $F\not\subseteq E_G(X)$, $d_{G'}(X)> d_G(X)-|F|$, and so the
first inequality in (\ref{e3.7}) is strict. Assume $F\subset E_G(X)$
below. If $|X|\geqslant 4$, then $E_G(X)$ is not a $\la''$-cut since
$G$ is super-$\la''$, which implies $d_G(X)
> \la''$, and so the second inequality in (\ref{e3.7}) is strict.

Now, consider $|X| = 3$ and we have the following two subcases.

If $g(G)\geqslant 4$, then $\la''=3k-4$, and so
$|F|=\la''-\xi+1=k-1\geqslant 3$. Since $F\subset E_G(X)$, there
exists one edge in $G[X]$ which is adjacent to at least two edges of
$F$, which implies $\xi(G')\leqslant \xi-2<\xi-1$, that is, the
inequality (\ref{e3.8}) is strict.

If $g(G)=3$, then $\la''=3k-6$. If $G[X]$ is not a triangle,
$d_G(X)=3k-4>3k-6=\la''$, and so the second inequality in (\ref{e3.7})
is strict. If $G[X]$ is a triangle, since
$|F|=\la''-\xi+1=k-3\geqslant 2$ and $F\subset E_G(X)$, then there
exists one edge in $G[X]$ which is adjacent to at least two edges of
$F$, which implies $\xi(G')\leqslant \xi-2<\xi-1$, that is, the
inequality (\ref{e3.8}) holds strictly.

Thus, the inequality (\ref{e3.6}) holds, and so  the inequality
(\ref{e3.5}) follows. We now prove the remaining parts of our
conclusions.

(i)\ When $g(G)\geqslant 4$, $\la''=3k-4$. By Theorem~\ref{thm1.6}
and (\ref{e3.5}), $k-1\geqslant \rho'(G)\geqslant \la''-\xi+1=k-1$,
which implies $\rho'(G)=k-1$.

(ii)\ When $g(G)=3$, $\la''=3k-6$. By (\ref{e3.5}),
$\rho'(G)\geqslant \la''-\xi+1=k-3$. If $k=5$, $\rho'(G)\geqslant
2$. Choose a subset $X\subset V(G)$ such that $G[X]$ is a triangle.
It is easy to check that $E_G(X)$ is a $\la''$-cut. Let $F$ be a set
of three independent edges of $E_G(X)$. Then $\la''(G-F)\leqslant
d_{G-F}(X)=6=\xi(G-F)$. This fact shows that $G-F$ is not
super-$\la'$, which implies $\rho'(G)\leqslant 2$. Thus,
$\rho'(G)=2$.

The theorem follows.
\end{pf}

\medskip

A graph $G$ is {\it transitive} if for any two given vertices $u$
and $v$ in $G$, there is an automorphism $\phi$ of $G$ such that
$\phi(u)=v$. A transitive graph is always regular. The studies on
extra edge-connected transitive graphs and super extra
edge-connected transitive graphs can be found in
\cite{m03,wl02,xx02,yzqg11} etc.

\begin{lemma}{\rm (Wang and Li~\cite{wl02})}\label{lem3.9}
Let $G$ be a connected transitive graph of degree $k\geqslant 4$
with girth $g\geqslant 5$. Then $G$ is $\la''$-optimal and
$\la''(G)=3k-4$.
\end{lemma}

\begin{lemma}{\rm (Yang {\it et al.}~\cite{yzqg11})}\label{lem3.10}
Let $G$ be a $C_4$-free transitive graph of degree $k\geqslant 4$.
If $G$ is $\la''$-optimal, then $G$ is super-$\la''$.
\end{lemma}

Combining Theorem~\ref{thm3.8} (i) with Lemma~\ref{lem3.9} and
Lemma~\ref{lem3.10}, we have the following corollary immediately.

\begin{coro}\label{cor3.11}
If $G$ is a connected transitive graph of degree $k\geqslant 4$ with
girth $g\geqslant 5$, then $\rho'(G)= k-1$.
\end{coro}

\begin{rem}\label{rem3.12}\ {\rm
In Corollary~\ref{cor3.11}, the condition ``\,$g\geqslant 5$" is
necessary. For example, the connected transitive graph $Q_4$ is
$\la''$-optimal and not super-$\la''$, and $\rho'(Q_4)=2$ (see
Remark~\ref{rem3.5}). }
\end{rem}

\section{$\bm{\rho'}$ for two families of networks}

As applications of Theorem~\ref{thm3.8} (i), in this section, we
determine the exact values of $\rho'(G)$ for two families of
networks $G(G_0,G_1;M)$ and $G(G_0,G_1,\dots,G_{m-1};\mathscr{M})$
subject to some conditions.

The first family of networks $G(G_0,G_1;M)$ is defined as follows. Let
$G_0$ and $G_1$ be two graphs with the same number of vertices. Then
$G(G_0,G_1;M)$ is the graph $G$ with vertex-set $V(G)=V(G_0)\cup
V(G_1)$ and edge-set $E(G)=E(G_0)\cup E(G_1)\cup M$, where $M$ is an
arbitrary perfect matching between vertices of $G_0$ and $G_1$. Thus
the hypercube $Q_n$, the twisted cube $TQ_n$, the crossed cube
$CQ_n$, the M\"{o}bius cube $MQ_n$ and the locally twisted cube
$LTQ_n$ all can be viewed as special cases of $G(G_0,G_1;M)$
(see~\cite{cth03}).

The second family of networks $G(G_0,G_1,\dots,G_{m-1};\mathscr{M})$ is
defined as follows. Let $G_0,G_1,\dots,G_{m-1}$ be $m~(\geqslant 3)$
graphs with the same number of vertices. Then $G(G_0,G_1,$
$\dots,G_{m-1};\mathscr{M})$ is the graph $G$ with vertex-set $V(G)
= V(G_0)\cup V(G_1)\cup \cdots\cup V(G_{m-1})$ and edge-set
$E(G)=E(G_0)\cup E(G_1)\cup \cdots\cup E(G_{m-1})\cup \mathscr{M}$,
where $\mathscr{M}=\cup_{i=0}^{m-1}M_{i,i+1({\rm mod}~m)}$ and
$M_{i,i+1({\rm mod}~m)}$ is an arbitrary perfect matching between
$V(G_i)$ and $V(G_{i+1({\rm mod}~m)})$. Recursive circulant graphs
\cite{pc94} and the undirected toroidal mesh \cite{x01} are special
cases of this family.

The super edge-connectivity of above two families of networks is
studied by Chen {\it et al.}~\cite{cth03}. Chen and Tan~\cite{ct07}
further studied the restricted edge-connectivity of above two
families of networks, and $\la'(G(G_0,G_1;M))$ is also studied by Xu
{\it et al.}~\cite{xww10}. The 2-extra edge-connectivity of above
two families of networks is studied by Wang {\it et
al.}~\cite{wyl08}. The vulnerability $\rho$ of super
edge-connectivity of the two families of networks is discussed by Wang
and Lu~\cite{wl12}. In this section, we will further investigate the
vulnerability $\rho'$ of the two families of super-$\la'$ networks
without triangles.

\begin{lemma}\label{lem4.1}{\rm (see Example 1.3.1 in Xu~\cite{x03})}
If $G$ is a triangle-free graph of order $n$, then $|E(G)|\leqslant
\frac{n^2}{4}$.
\end{lemma}

We consider the first family of graphs $G=G(G_0,G_1;M)$ for
$k$-regular triangle-free and super-$\la$ graphs $G_0$ and $G_1$.
Under these hypothesis, $G$ is $(k+1)$-regular and triangle-free. By
Theorem~\ref{thm3.2}, we can assume $k\geqslant 3$. We attempt to
use Theorem~\ref{thm3.8} (i) to determine the exact value of
$\rho'(G)$ when $G$ is super-$\la''$. However, there are some such
graphs that are not super-$\la''$.

\begin{exam}{\rm
Let $G_0$ be a $k$-regular triangle-free and super-$\la$ graph of
order $n$. Then $G_0$ is $\la'$-connected, $k\geqslant 3$ and $n\geqslant 6$.
$G=G_0\times K_2$ can be viewed as $G(G_0,G_0;M)$ for some perfect matching $M$.
Assume $n\leqslant 3k-1$ or $\la'(G_0)\leqslant \frac{3k-1}{2}$. If
the former happens, then $M$ is a $2$-extra edge-cut, and so
$\la''(G_0\times K_2)\leqslant |M|=3k-1$. However, $G_0\times K_2$
is not super-$\la''$ since $n\geqslant 6$. If the latter happens,
let $X_0\subset V(G_0)$ such that $E_{G_0}(X_0)$ is a $\la'$-cut of
$G_0$, then $G[X_0]\times K_2\subset G$. Let $Y=V(G[X_0]\times
K_2)$. Since $E_G(Y)$ is a $2$-extra edge-cut, $\la''(G_0\times
K_2)\leqslant |E_G(Y)|= 2\la'(G_0)\leqslant 3k-1$. However,
$G_0\times K_2$ is not super-$\la''$ since $|Y|\geqslant 4$. }
\end{exam}

This example shows that the condition
``\,$\min\{n,\la'(G_0)+\la'(G_1)\}>3k-1$\," is necessary to
guarantee that $G=G(G_0,G_1;M)$ is super-$\la''$. Thus, we can state
our result as follows.

\begin{theorem}\label{thm4.3}
Let $G_i$ be a triangle-free $k$-regular and super-$\la$ graph of
order $n$ for each $i=0,1$. If $\min\{n,\la'_0+\la'_1\}>3k-1$, then
$G=G(G_0,G_1;M)$ is super-$\la''$ and $\rho'(G)=k$, where
$\la'_i=\la'(G_i)$ for each $i=0,1$.
\end{theorem}

\begin{pf}
Clearly, $k\geqslant 3$. Since $G$ is $(k+1)$-regular and
triangle-free, by Theorem~\ref{thm3.8} (i), we only need to prove
that $G$ is super-$\la''$. Since $M$ is a 2-extra edge-cut of $G$,
$\la''(G)$ exists. By Theorem \ref{thm1.2},
\begin{equation}\label{e4.1}
 \la''(G)\leqslant \xi_2(G)=3k-1.
 \end{equation}
Suppose to the contrary that $G$ is not super-$\la''$. Then there
exists a $\la''$-fragment $X$ of $G$ such that
$|\overline{X}|\geqslant |X|\geqslant 4$. Since $G$ is
triangle-free, $G[X]$ is also triangle-free, and so
$|E(G[X])|\leqslant \frac{|X|^2}{4}$ by Lemma \ref{lem4.1}. It
follows that
 $$
 \begin{array}{rl}
 3k-1\geqslant d_G(X)=(k+1)|X|-2|E(G[X])|\geqslant (k+1)|X|-\frac{1}{2}|X|^2,
 \end{array}
 $$
that is, $ (|X|-3)(|X|-(2k-1))+1\geqslant 0$, which implies that,
since $|X|\geqslant 4$ and $k\geqslant 3$,
 \begin{equation}\label{e4.2}
 |X|\geqslant 2k-1.
 \end{equation}

We will deduce a contradiction to (\ref{e4.1}) by proving that
 \begin{equation}\label{e4.3}
  \la''(G)>3k-1.
 \end{equation}
To the end, set $V_i=V(G_i)$ and $X_i=X\cap V_i$ for each $i=0,1$.
There are two cases.

\vskip6pt

{\it Case 1.} Exactly one of $X_0$ and $X_1$ is empty.

Without loss of generality, assume $X=X_0$. Then
$E_G(X)=E_{G_0}(X)\cup [X,V_1]$. By the definition of $G$,
$|[X,V_1]|=|X|$, and so
 \begin{align}\label{e4.4}
 \la''(G)=d_G(X)=d_{G_0}(X)+|X|.
 \end{align}
It is easy to check that $G_0[V_0\setminus X]$ is connected. Thus,
when $2\leqslant |X|\leqslant n-2$, $E_{G_0}(X)$ is a 1-extra
edge-cut of $G_0$, and so $d_{G_0}(X)\geqslant \la'_0$. Since $G_0$
is super-$\la$, $\la'_0>\la(G_0)=k$, and so
 \begin{align}\label{e4.5}
 d_{G_0}(X)\geqslant \left\{
 \begin{array}{ll}
k+1 &\ {\rm if}\ 2\leqslant |X|\leqslant n-2,\\
k &\ {\rm if}\ |X|=n-1,\\
0            &\ {\rm if}\ |X|=n.
\end{array}
\right.
\end{align}
Substituting $n>3k-1$, (\ref{e4.2}) and (\ref{e4.5}) into
(\ref{e4.4}) yields the inequality (\ref{e4.3}).

\vskip6pt

{\it Case 2.} $X_0\ne \emptyset$ and $X_1\ne \emptyset$.

Assume that one of $G[X_0]$, $G[X_1]$, $G[V_0\setminus X_0]$ and
$G[V_1\setminus X_1]$ is not connected. Without loss of generality,
assume that $G[X_0]$ has two components $H$ and $T$. Then
$[H,V_0\setminus X_0]\cup [T,V_0\setminus X_0]\cup [X_1,V_1\setminus
X_1]\subseteq E_G(X)$, and the first two are edge-cuts of $G_0$, and
the last is an edge-cut of $G_1$. Since $G_i$ is super-$\la$,
$\la(G_i)=k$ for each $i=0,1$. Thus,
 $$
 \la''(G)=|E_G(X)|\geqslant |[H,V_0\setminus X_0]|+|[T,V_0\setminus X_0]|+|[X_1,V_1\setminus
 X_1]|\geqslant 3k>3k-1,
 $$
and so (\ref{e4.3}) follows.

Now, we assume that all of $G[X_0]$, $G[X_1]$, $G[V_0\setminus X_0]$
and $G[V_1\setminus X_1]$ are connected. Since $|X|\geqslant 4$,
$\max\{|X_0|,|X_1|\}\geqslant 2$. We consider the following two
subcases.

\vskip6pt

{\it Subcase $2.1$.} $|X_0|\geqslant 2$ and $|X_1|\geqslant 2$.

In this case, $E_{G_0}(X_0)\cup E_{G_1}(X_1)\subseteq E_G(X)$. For
each $i=0,1$, $d_{G_i}(X_i)\geqslant \la'_i$ since $E_{G_i}(X_i)$ is
a $1$-extra edge-cut of $G_i$. By our hypothesis,
 $$
 \la''(G)=d_G(X)\geqslant d_{G_0}(X_0)+d_{G_1}(X_1)\geqslant \la'_0+\la'_1>3k-1,
 $$
and so (\ref{e4.3}) follows.

\vskip6pt

{\it Subcase $2.2$.} Exact one of $X_0$ and $X_1$ is a single
vertex.

Without loss of generality, assume $|X_0|=1$. Then
$|X_1|=|X|-1\geqslant 2k-2$ by (\ref{e4.2}). Clearly,
 $$
 E_G(X)=E_{G_0}(X_0)\cup E_{G_1}(X_1)\cup [X_1,V_0\setminus X_0],
 $$
$d_{G_0}(X_0)=k$ and $|[X_1,V_0\setminus X_0]|=|X|-2\geqslant
2k-3$, and so
  \begin{align}\label{e4.6}
 \la''(G)=d_G(X)\geqslant k+d_{G_1}(X_1)+2k-3.
  \end{align}
If $2\leqslant |X_1|\leqslant n-2$, then $X_1$ is a 1-extra edge-cut
of $G_1$, and so $d_{G_1}(X_1)\geqslant \la'_1>\la(G_1)=k$ since
$G_1$ is super-$\la$. If $|X_1|=n-1$, then $E_{G_1}(X_1)$ isolates a
vertex, and so $d_{G_1}(X_1)=k$. Thus, we always have
$d_{G_1}(X_1)\geqslant k$. Substituting this inequality into
(\ref{e4.6}) yields (\ref{e4.3}) since $d_G(X)\geqslant
k+k+2k-3=4k-3>3k-1$ for $k\geqslant 3$.

Under the hypothesis that $G$ is not super-$\la''$, we deduce a
contradiction to (\ref{e4.1}). Thus, $G$ is super-$\la''$. By
Theorem~\ref{thm3.8} (i), $\rho'(G)=k$, and so the theorem follows.
\end{pf}

\begin{lemma}{\rm (Xu {\it et al.}~\cite{xww10})}\label{lem4.4}
If $G_n\in \{Q_n, TQ_n, CQ_n, MQ_n, LTQ_n\}$, then  $\la'(G_n)=2n-2$
and, thus, $G_n$ is $\la'$-optimal for $n\geqslant 2$, and is
super-$\la$ for $n\geqslant 3$.
\end{lemma}

\begin{coro}\label{cor4.5}
Let $G_n\in \{Q_n, TQ_n, CQ_n, MQ_n, LTQ_n\}$. If $n\geqslant 5$,
then $G_n$ is super-$\la''$, super-$\la'$ and $\rho'(G_n)=n-1$.
\end{coro}

\begin{pf}
Let $G_n\in \{Q_n, TQ_n, CQ_n, MQ_n, LTQ_n\}$. Then $G_n$ can be
viewed as the graph $G(G_{n-1},G_{n-1};M)$ corresponding to some
perfect matching $M$. $G_{n-1}$ is an $(n-1)$-regular and
triangle-free graph of order $2^{n-1}$. By Lemma \ref{lem4.4},
$G_{n-1}$ is super-$\la$ and $\la'(G_{n-1})=2n-4$ for $n\geqslant
4$. Thus, $2\la'(G_{n-1})=4n-8>3(n-1)-1$ and $2^{n-1}>3(n-1)-1$ for
$n\geqslant 5$. By Theorem \ref{thm4.3}, $G_n$ is super-$\la''$ and
$\rho'(G_n)=n-1$ if $n\geqslant 5$. Hence, if $n\geqslant 5$,
$\la''(G_n)=3n-4>2n-2=\xi(G_n)$ implies $G_n$ is super-$\la'$.
\end{pf}

\vskip6pt

\begin{rem}\label{rem4.6}\ {\rm
In Corollary~\ref{cor4.5}, the condition ``\,$n\geqslant 5$\," is
necessary. For example, $Q_4$ is $\la''$-optimal and not
super-$\la''$, and $\rho'(Q_4)=2$ (see Remark~\ref{rem3.5}). }
\end{rem}

\vskip6pt

We now consider the second family of graphs
$G(G_0,G_1,\dots,G_{m-1};\mathscr{M})$. To guarantee that $G$ is
triangle-free, we can assume $m\geqslant 4$. Let
$I_m=\{0,1,\dots,m-1\}$.

\begin{theorem}\label{thm4.7}
Let $G_i$ be a $k$-regular $k$-edge-connected graph of order $n$
without triangles for each $i\in I_m$. If $k\geqslant 3$, $n>
\lceil\frac{3k+2}{2}\rceil$ and $m\geqslant 4$, then
$G=G(G_0,\dots,G_{m-1};\mathscr{M})$ is super-$\la''$ and
$\rho'(G)=k+1$.
\end{theorem}

\begin{pf}
It is easy to check that $G$ is $(k+2)$-regular and triangle-free.
By Theorem \ref{thm1.2}, $G$ is $\la''$-connected and
\begin{equation}\label{e4.7}
 \la''(G)\leqslant \xi_2(G)=3k+2.
 \end{equation}
By Theorem~\ref{thm3.8} (i), we only need to prove that $G$ is
super-$\la''$. Suppose to the contrary that $G$ is not
super-$\la''$. Then there exists a $\la''$-fragment $X$ of $G$ such
that $|\overline X|\geqslant |X|\geqslant 4$. Since $G$ is
triangle-free, $G[X]$ is also triangle-free and $|E(G[X])|\leqslant
\frac{|X|^2}{4}$ by Lemma \ref{lem4.1}. It follows that
 $$
 \begin{array}{rl}
 3k+2\geqslant \la''(G)=d_G(X)=(k+2)|X|-2|E(G[X])|\geqslant (k+2)|X|-\frac{1}{2}|X|^2,
 \end{array}
 $$
that is, $(|X|-3)(|X|-(2k+1))+1\geqslant 0$, which implies, since
$|X|\geqslant 4$ and $k\geqslant 3$,
 \begin{equation}\label{e4.8}
 |X|\geqslant 2k+1.
 \end{equation}

We will deduce a contradiction to (\ref{e4.7}) by proving that
\begin{equation}\label{e4.9}
 \la''(G)>3k+2.
 \end{equation}

To the end, for each $i\in I_m$, let
 $$
 \begin{array}{l}
 V_i=V(G_i),\ \ X_i=X\cap V_i,\\
 F_i=E_G(X)\cap E(G_i),\ \
 F'_i=E_G(X)\cap M_{i,i+1({\rm mod}~m)}.
 \end{array}
 $$
Then $F_i=E_{G_i}(X_i)$.  Let
 $$
 J=\{j\in I_m:\ X_j\ne \emptyset\}\ {\rm and}\ J'=\{j\in J:\ X_j=V_j\}.
 $$
Then $|F_j|\geqslant \la(G_j)=k$ for any $j\in J\setminus J'$. Thus,
if $|J\setminus J'|\geqslant 4$, then
 $$
  \la''(G)=|E_G(X)|\geqslant \sum_{i\in J\setminus J'}|F_i|\geqslant 4k> 3k+2\ \ {\rm for}\ k\geqslant 3,
 $$
and so (\ref{e4.9}) follows. We assume $|J\setminus J'|\leqslant 3$
below. There are two cases.

\vskip6pt

{\it Case} 1. 
$|J|\leqslant m-1$.

\vskip6pt

{\it Subcase} 1.1. $J'\ne \emptyset$.

Let $\ell\in I_m\setminus J$ and $j\in J'$. Then $X_\ell=\emptyset$
and $X_j=V_j$. Since $j\ne\ell$, without loss of generality, assume
$\ell<j$ and let $s=j-\ell$. By the structure of $G$, there exist
exactly $n$ disjoint paths of length $s$ between $V_j$ and $V_\ell$
passing through $G_{j-1}$ (maybe $j-1=\ell$), and $n$ disjoint paths
of length $m-s$ between $V_j$ and $V_\ell$ passing through
$G_{j+1({\rm mod}~m)}$ (maybe $\ell=j+1({\rm mod}~m)$). Each of
these paths has at least one edge that is in $E_G(X)$. Since $n>
\lceil\frac{3k+2}{2}\rceil$, we have that
 $$
 \la''(G)=|E_G(X)|\geqslant 2n>3k+2,
 $$
and so (\ref{e4.9}) follows.

\vskip6pt

{\it Subcase} 1.2. $J'=\emptyset$. In this subcase, $|J|\leqslant
3$.

If $|J|=1$, say $X_1=X$, since $E_G(X)=F_1\cup F'_0\cup F'_1$ and
$|F'_0|=|F'_1|=|X_1|=|X|$, then $|E_G(X)|\geqslant 2|X|+|F_1|$.
Combining this fact with (\ref{e4.8}), we have that
 $$
 \la''(G)=|E_G(X)|\geqslant 2|X|+|F_1|\geqslant 2(2k+1)+k>3k+2,
 $$
and so (\ref{e4.9}) follows.

If $|J|=2$, say $J=\{p,q\}$, then $|p-q|=1$ since $G[X]$ is
connected, say $q=p+1$. $F_p\cup F_{p+1}\cup F'_{p-1}\cup
F'_{p+1}\subseteq E_G(X)$. Since $|F_p|\geqslant k$,
$|F_{p+1}|\geqslant k$ and $|F'_{p-1}\cup F'_{p+1}|=|X|$. Combining
these facts with (\ref{e4.8}), we have that
 $$
 \la''(G)=|E_G(X)|\geqslant |F_p|+|F_{p+1}|+|X|\geqslant 2k+(2k+1)>3k+2,
 $$
and so (\ref{e4.9}) follows.

If $|J|=3$, without loss of generality, assume $J=\{1,2,3\}$ since
$G[X]$ is connected, then $F'_0\ne\emptyset$ and $F'_3\ne\emptyset$
since $m\geqslant 4$. If $|F'_0|=|F'_3|=1$, then $|X_1|=|X_3|=1$.
Since $|X|\geqslant 4$, if $|X_1|=|X_3|=1$, then $|X_2|\geqslant 2$,
and so $|F'_1|\geqslant 1$ and $|F'_2|\geqslant 1$. Thus, it is
always true that $|F'_0|+|F'_1|+|F'_2|+|F'_3|>2$. It follows that
 $$
 \la''(G)=|E_G(X)|\geqslant \sum_{j=1}^3|F_j|+\sum_{j=0}^3|F'_j|>
  3k+2,
 $$
and so (\ref{e4.9}) follows.

\vskip6pt

 {\it Case} 2. $|J|=m$.

In this case, $J'\ne \emptyset$ since $m\geqslant 4$ and
$|J\setminus J'|\leqslant 3$. If $|J\setminus J'|\leqslant 2$, then
$|J'|\geqslant 2$, and
 $$
 |\overline X|=\sum\limits_{j\in J\setminus J'}|V(G_j-X_j)|\leqslant
 2(n-1)<2n\leqslant \sum\limits_{j\in J'}|V_j|<|X|,
 $$
a contradiction to $|\overline X|\geqslant |X|$. Therefore,
$|J\setminus J'|=3$. Since $G[X]$ is connected, without loss of
generality, let $J\setminus J'=\{1,2,3\}$. Then $F'_0\ne\emptyset$
and $F'_3\ne\emptyset$. Since $|\overline X|\geqslant |X|$, there
exists at least two $i\in J\setminus J'$ such that $|X_i|\leqslant
\frac{n}{2}$. Thus, at least one of $|V_1\setminus X_1|$ and
$|V_3\setminus X_3|$ is not less than $\frac{n}{2}$, that is, either
$|F'_0|\geqslant \frac{n}{2}$ or $|F'_3|\geqslant \frac{n}{2}$. It
follows that
 $$
  \la''(G)=|E_G(X)|\geqslant \sum_{j=1}^3|F_j|+|F'_0|+|F'_3|> 3k+2,
 $$
and so (\ref{e4.9}) follows.

Under the hypothesis that $G$ is not super-$\la''$, we deduce a
contradiction to (\ref{e4.7}). Thus, $G$ is super-$\la''$. By
Theorem~\ref{thm3.8} (i), $\rho'(G)=k+1$, and so the theorem
follows.
\end{pf}

\bigskip

As applications of Theorem \ref{thm4.7}, we consider two families
well-known transitive networks.

Let $G(n,d)$ denote a graph which has the vertex-set
$V=\{0,1,\dots,n-1\}$, and two vertices $u$ and $v$ are adjacent if
and only if $|u-v|=d^i({\rm mod}~n)$ for any $i$ with $0\leqslant
i\leqslant \lceil\log_d n\rceil-1$. Clearly, $G(d^m,d)$ is a
circulant graph, which is $\delta$-regular and $\delta$-connected,
where $\delta=2m-1$ if $d=2$ and $\delta=2m$ if $d\ne 2$. For
circulant graphs with order between $d^{m}$ and $d^{m+1}$, that is,
$G(cd^{m},d)$ with $1<c<d$, $\delta=2m+1$ if $c=2$ and $\delta=2m+2$
if $c>2$, moreover, Park and Chwa~\cite{pc94} showed that
$G(cd^{m},d)$ can be recursively constructed, that is,
$G(cd^{m},d)=G(G_0,G_1,\dots,G_{d-1};\mathscr{M})$, where $G_i$ is
isomorphic to $G(cd^{m-1},d)$ for each $i=0,1,\ldots,d-1$, and so
$G(cd^{m},d)$ is called the {\it recursive circulant graph}, which
is $\delta$-regular and $\delta$-connected. In particular, the graph
$G(2^m,4)$ is $2m$-regular $2m$-connected, has the same number of
vertices and edges as a hypercube $Q_m$. However, $G(2^m,4)$ with
$m\geqslant 3$ is not isomorphic to $Q_m$ since $G(2^m,4)$ has an
odd cycle of length larger than $3$. Compared with $Q_m$, $G(2^m,
4)$ achieves noticeable improvements in diameter
($\lceil\frac{3m-1}{4}\rceil$). Thus, the recursive circulant graphs
have attracted much research interest in recent ten years (see
Park~\cite{p08,pc94} and references therein). Since, when
$c\geqslant 3$, $G(cd^0,d)$ is isomorphic to a cycle of length $c$,
$G(cd^m,d)$ contains triangles if $c=3$.

The {\it $n$-dimensional undirected toroidal mesh}, denoted by
$C(d_1,\dots,d_n)$, is defined as the cartesian products
$C_{d_1}\times C_{d_2}\times\cdots\times C_{d_n}$, where $C_{d_i}$
is a cycle of length $d_i~(\geqslant 3)$ for each $i = 1, 2,\dots,
n$ and $n\geqslant 2$. It is known that $C(d_1,\dots,d_n)$ is a
$2n$-regular $2n$-edge-connected transitive graph with girth $g =
\min\{4, d_i, 1 \leqslant i \leqslant n\}$. Thus, if $d_i\geqslant
4$ for each $i=1,2,\ldots, n$, then $C(d_1,\dots,d_n)$ is
triangle-free.  $C(d_1,\dots,d_n)$ can be expressed as
$G(G_0,G_1,\dots,G_{d_1-1};\mathscr{M})$, where $G_i$ is isomorphic
to $C_{d_2}\times\cdots\times C_{d_n}$ for each
$i=0,1,\ldots,d_1-1$.

Since the two families of networks are transitive, by
Corollary~\ref{cor3.11}, we can determine the exact values of $\rho'$
when the girth $g\geqslant 5$. By Theorem \ref{thm4.7}, we have the
following two stronger results immediately.

\begin{coro}
Let $c,d,m$ be three positive integers with $1<c<d$, $c\ne 3$,
$d\geqslant 4$, $m\geqslant 2$. Then $G=G(cd^{m},d)$ is super-$\la''$,
and $\rho'(G)=2m$ if $c=2$ and $\rho'(G)=2m+1$ if $c\geqslant 4$.
\end{coro}

\begin{coro}
If $n\geqslant 3$ and $d_i\geqslant 4$ for $1 \leqslant i \leqslant
n$, then $G=C(d_1,\dots,d_n)$ is super-$\la''$ and
$\rho'(G)=2n-1$.
\end{coro}

\end{document}